\documentclass[12pt,draft]{article}
\usepackage{amssymb,amsmath,amsthm,latexsym}
\usepackage[cp1251]{inputenc}
\usepackage[russian]{babel}

\begin{document}

\begin{center}{\large \bf Описание алгоритма построения интегрального представления по вееру торического многообразия}\end{center}
\begin{center}
{\hfill\large\bf А.А.Кытманов\footnote{${}^*$ \copyright\ А.А.Кытманов, Сибирский федеральный университет, 2010.}}
\end{center}

\begin{abstract}
В настоящей работе приводится описание алгоритма компьютерной алгебры, реализующего построение торического многообразия по его вееру и построение интегрального представления в пространстве $\mathbb{C}^d$, связанного с данным торическим многообразием. Алгоритмы были разработаны с помощью методов, описанных в работе \cite{AAK}.
\end{abstract}

Опишем процедуру {\bf ToricVariety}, реализующую конструкцию торического многообразия и связанного с ним интегрального представления. Входными данными процедуры являются набор векторов (одномерных образующих) веера $\Sigma$ и набор из номеров векторов, образующих конусы максимальной размерности в $\Sigma$.

\medskip

\noindent {\bf ToricVariety} := {\bf procedure}($vec\_list$, $cone\_list$)

\noindent \hskip3mm d := число элементов в $vec\_list$

\noindent \hskip3mm n := число элементов в первом элементе из $vec\_list$

\noindent \hskip3mm $V$ := пустой список, $Z$ := пустой список

\noindent \hskip3mm $Z2$ := пустой список, $Dz$ := пустой список

\noindent \hskip3mm {\bf for} $i$ {\bf from} $1$ {\bf to} $d$ {\bf do}

\noindent \hskip6mm добавить в $V$ элемент $v[i]$, добавить в $Z$ элемент $z[i]$, 

\noindent \hskip6mm добавить в $Z2$ элемент $|z[i]|^2$, добавить в $Dz$ элемент $dz[i]$

\noindent \hskip3mm {\bf end do}

\noindent \hskip3mm {\bf print}(``Рассмотрим торическое многообразие, веер которого задан векторами'')

\noindent \hskip3mm {\bf for} $i$ {\bf from} 1 {\bf to} $d$ {\bf do}

\noindent \hskip6mm {\bf print}($v[i]=vec\_list[i]$)

\noindent \hskip3mm {\bf end do}

\noindent \hskip3mm $CL$ := пустой список

\noindent \hskip3mm {\bf for} $i$ {\bf from} 1 {\bf to} число элементов в $cone\_list$ {\bf do}

\noindent \hskip6mm $Cl$ := пустой список

\noindent \hskip6mm {\bf for} $j$ {\bf from} 1 {\bf to} число элементов в $cone\_list[i]$ {\bf do}

\noindent \hskip9mm $t[j]$ := $v[cone\_list[i][j]]$

\noindent \hskip9mm добавить в $Cl$ элемент t[j]

\noindent \hskip6mm {\bf end do}

\noindent \hskip6mm добавить в $CL$ элемент $Cl$

\noindent \hskip3mm {\bf end do}

\noindent \hskip3mm {\bf print}(``Конусы максимальной размерности порождены векторами:'')

\noindent \hskip3mm {\bf print}(элементы $CL$)

\medskip
\noindent \% Конструкция торического многообразия
\medskip

\noindent \hskip3mm LR := {\bf LinRel}($vec\_list$)

\noindent \hskip3mm $V\_LR$ := пустой список, $W$ := пустой список

\noindent \hskip3mm {\bf for} $i$ {\bf from} 1 {\bf to} число элементов в $LR$ {\bf do}

\noindent \hskip6mm добавить в $V\_LR$ cумму элементов списка $LR[i]$, умноженных на соответствующие элементы списка $V$ 

\noindent \hskip6mm добавить в $W$ cумму элементов списка $LR[i]$, умноженных на соответствующие элементы списка $Z2$

\noindent \hskip3mm {\bf end do}

\noindent \hskip3mm $P$ := пустой список

\noindent \hskip3mm {\bf for} $i$ {\bf from} 1 {\bf to} число элементов в $LR$ {\bf do}

\noindent \hskip6mm добавить в $P$ элемент $\rho[i]$

\noindent \hskip3mm {\bf end do}

\noindent \hskip3mm $\rho$ := вектор из элементов списка $P$

\noindent \hskip3mm {\bf print}(``Все независимые линейные соотношения на вектора заданного веера задаются системой:'')

\noindent \hskip3mm {\bf for} $i$ {\bf from} 1 {\bf to} число элементов в $LR$ {\bf do}

\noindent \hskip6mm {\bf print}($V\_LR[i]=0$)

\noindent \hskip3mm {\bf end do}

\noindent \hskip3mm $M$ := {\bf Matrix}($LR$)

\noindent \hskip3mm $vv$ := общее решение однородной СЛАУ, задающейся матрицей $M$ со свободными переменными $v_i$

\noindent \hskip3mm $VV$ := пустой список

\noindent \hskip3mm {\bf for} $i$ {\bf from} 1 {\bf to} $d$ {\bf do}

\noindent \hskip6mm добавить в $VV$ элемент $v[i]=vv[i]$

\noindent \hskip3mm {\bf end do}

\noindent \hskip3mm $L\_power$ := список из столбцов матрицы $M$

\noindent \hskip3mm $\Lambda$ := пустой список

\noindent \hskip3mm {\bf for} $i$ {\bf from} 1 {\bf to} $d-n$ {\bf do}

\noindent \hskip6mm добавить в $\Lambda$ элемент $\lambda[i]$

\noindent \hskip3mm {\bf end do}

\noindent \hskip3mm $G$ := пустой список

\noindent \hskip3mm {\bf for} $i$ {\bf from} 1 {\bf to} $d$ {\bf do}

\noindent \hskip6mm добавить в $G$ произведение элементов списка $\Lambda$, возведенных в степени соответствующих элементов списка $L\_power[i]$

\noindent \hskip3mm {\bf end do}

\noindent \hskip3mm $PC$ := {\bf PrimColl}($cone\_list$)

\noindent \hskip3mm {\bf print}(``Заданное торическое многообразие $\mathbb{X}$ является фактор-пространством:'')

\noindent \hskip3mm {\bf print}(X=($\mathbb{C}^d\setminus Z(\Sigma))/G$)

\noindent \hskip3mm {\bf print}(``Здесь сингулярное множество $Z(\Sigma)$ есть объединение координатных подпространств $z[i[1]]=\ldots=z[i[n]]=0$, где $\{i[1],\ldots,i[n]\}\in PC$'')

\noindent \hskip3mm {\bf print}(``Группа $G$ является $d-n$-параметрической поверхностью в $\mathbb{C}^d$'')

\noindent \hskip3mm {\bf print}($G$)

\medskip
\noindent \% Конструкция конуса Кэлера
\medskip

\noindent \hskip3mm $Z\_s$ := пустой список, $Z\_v$ := пустой список

\noindent \hskip3mm {\bf for} i {\bf from} 1 {\bf to} число элементов в $PC$ {\bf do}

\noindent \hskip6mm $z\_s$ := пустой список, $z\_v$ := пустой список

\noindent \hskip6mm {\bf for} j {\bf from} 1 {\bf to} число элементов в $PC[i]$ {\bf do}

\noindent \hskip9mm $ts[j]$ := $|z[PC[i][j]]|^2$

\noindent \hskip9mm $t[j]$ := $v[PC[i][j]]$

\noindent \hskip9mm добавить в $z\_s$ элемент $ts[j]$

\noindent \hskip9mm добавить в $z\_v$ элемент $t[j]$

\noindent \hskip6mm {\bf end do}

\noindent \hskip6mm добавить в $Z\_s$ элемент $z\_s$

\noindent \hskip6mm добавить в $Z\_v$ элемент $z\_v$

\noindent \hskip3mm {\bf end do}

\noindent \hskip3mm $z\_s$ := пустой список, $z\_v$ := пустой список

\noindent \hskip3mm {\bf for} $i$ {\bf from} 1 {\bf to} число элементов в $Z\_v$ {\bf do}

\noindent \hskip6mm добавить в $z\_s$ сумму элементов списка $Z\_s[i]$

\noindent \hskip6mm добавить в $z\_v$ сумму элементов списка $Z\_v[i]$

\noindent \hskip3mm {\bf end do}

\noindent \hskip3mm {\bf for} $i$ {\bf from} 1 {\bf to} $d$ {\bf do}

\noindent \hskip6mm $v[i]$ := $vec\_list[i]$

\noindent \hskip3mm {\bf end do}

\noindent \hskip3mm $zero$ := пустой список

\noindent \hskip3mm {\bf for} $i$ {\bf from} 1 {\bf to} $n$ {\bf do}

\noindent \hskip6mm добавить в $zero$ элемент 0

\noindent \hskip3mm {\bf end do}

\noindent \hskip3mm $Kahler$ := пустой список

\noindent \hskip3mm {\bf for} $i$ {\bf from} 1 {\bf to} число элементов в $z\_v$ {\bf do}

\noindent \hskip6mm $ss$ := 0

\noindent \hskip6mm {\bf if} (значение $z\_v[i]$)=zero {\bf then}

\noindent \hskip9mm {\bf for} j {\bf from} 1 {\bf to} число элементов в $LR$ {\bf do}

\noindent \hskip12mm {\bf if} $z\_v[i]=V\_LR[j]$ {\bf then} добавить в $Kahler$ элемент $rho[j]$ {\bf end if}

\noindent \hskip9mm {\bf end do}

\noindent \hskip6mm {\bf else}

\noindent \hskip9mm {\bf for} $j$ {\bf from} 1 {\bf to} число элементов в $CL$ {\bf do}

\noindent \hskip12mm $t$ := элементы ({\bf LinRel}([(-1)*(значение $z\_v[i]$),(значение $CL[j]$)]))

\noindent \hskip12mm $s$ := 0

\noindent \hskip12mm {\bf for} $k$ {\bf from} 1 {\bf to} число элементов в $t$ {\bf do}

\noindent \hskip15mm {\bf if} $t[k]<0$ {\bf then} $s$ := 1 {\bf end if}

\noindent \hskip12mm {\bf end do}

\noindent \hskip12mm {\bf if} $s=0$ {\bf and} $ss=0$ {\bf then}

\noindent \hskip15mm $tt$ := [$z\_v[i]$, $(-1)$*$CL[j]$]

\noindent \hskip15mm $s$ := сумма элементов списка $t$, умноженных на соответствующие элементы списка $tt$ 

\noindent \hskip15mm добавить в $Kahler$ элементы подстановки $VV$ в выражение $s$ 

\noindent \hskip15mm $ss$ := 1 

\noindent \hskip12mm {\bf end if}

\noindent \hskip9mm {\bf end do}

\noindent \hskip6mm {\bf end if}

\noindent \hskip3mm {\bf end do}

\noindent \hskip3mm {\bf print}(``Конус Кэлера данного торического многообразия определяется следующей системой неравенств:'')

\noindent \hskip3mm {\bf for} $i$ {\bf from} 1 {\bf to} число элементов в $Kahler$ {\bf do}

\noindent \hskip6mm {\bf print}($Kahler[i]>0$)

\noindent \hskip3mm {\bf end do}

\medskip
\noindent \% Вычисление формы $h(z)$ в числителе дифференциальной формы $\omega$
\medskip

\noindent \hskip3mm $H$ := 0

\noindent \hskip3mm $Dd$ := пустое множество

\noindent \hskip3mm {\bf for} $i$ {\bf from} 1 {\bf to} $d$ {\bf do}

\noindent \hskip6mm добавить в $Dd$ элемент $i$

\noindent \hskip3mm {\bf end do}

\noindent \hskip3mm $s$ := {\bf SetOfAllSubsets}($Dd$, $n$)

\noindent \hskip3mm {\bf for} $i$ {\bf from} 1 {\bf to} число элементов в $s$ {\bf do}

\noindent \hskip6mm $ss$ := $(-1)^{(\text{сумма элементов множества}\ s[i]-1)}$

\noindent \hskip6mm $t$ := разность множеств $Dd$ и $s[i]$

\noindent \hskip6mm $sz$ := 1

\noindent \hskip6mm {\bf for} $j$ {\bf from} 1 {\bf to} число элементов в $s[i]$ {\bf do}

\noindent \hskip9mm $sz$ := $sz\wedge dz[s[i][j]]$

\noindent \hskip6mm {\bf end do}

\noindent \hskip6mm $tt$ := пустой список

\noindent \hskip6mm $tz$ := 1

\noindent \hskip6mm {\bf for} $j$ {\bf from} 1 {\bf to} число элементов в $t$ {\bf do}

\noindent \hskip9mm добавить в $tt$ элемент $L\_power[t[j]]$

\noindent \hskip9mm $tz$ := $tz*z[t[j]]$

\noindent \hskip6mm {\bf end do}

\noindent \hskip6mm $A$ := определитель матрицы, построенной на векторах из списка $tt$

\noindent \hskip6mm $H$ := $H+ss*A*tz*sz$

\noindent \hskip3mm {\bf end do}

\medskip
\noindent \% Вычисление знаменателя дифференциальной формы $\omega$
\medskip

\noindent \hskip3mm $tt$ := пустой список

\noindent \hskip3mm $ss$ := пустой список

\noindent \hskip3mm $g$ := 0

\noindent \hskip3mm {\bf for} $i$ {\bf from} 1 {\bf to} число элементов в $cone\_list$ {\bf do}

\noindent \hskip6mm добавить в $tt$ разность множеств $Dd$ и $cone\_list[i]$

\noindent \hskip6mm $s$ := пустое множество

\noindent \hskip6mm $t$ := 1

\noindent \hskip6mm {\bf for} $j$ {\bf from} 1 {\bf to} число элементов в $tt[i]$ {\bf do}

\noindent \hskip9mm добавить в $s$ элемент $|z[tt[i][j]]|^{2*({\bf nu\_sigma}(vec\_list,Permut(cone\_list[i],tt[i][j]))+1)}$

\noindent \hskip9mm $t$ := $t*s[j]$

\noindent \hskip6mm {\bf end do}

\noindent \hskip6mm добавить в $ss$ элемент $s$

\noindent \hskip6mm $g$ := $g+t$

\noindent \hskip3mm {\bf end do}

\noindent \hskip3mm {\bf print}(``Ядром интегрального представления является дифференциальная форма'')

\noindent \hskip3mm {\bf print}($\omega(z)=\overline{h(z)}\wedge dz/g(z,\overline{z})$)

\noindent \hskip3mm {\bf print}(`` где форма $h(z)=$'')

\noindent \hskip3mm {\bf print}($H$)

\noindent \hskip3mm {\bf print}(`` а $dz$ --- форма:'')

\noindent \hskip3mm {\bf print}(внешнее произведение элементов $Dz$)

\noindent \hskip3mm {\bf print}(`` Знаменатель $g(z,\overline{z})$ --- функция'')

\noindent \hskip3mm {\bf print}($g$)

\noindent \hskip3mm {\bf print}(``ОСНОВНОЙ РЕЗУЛЬТАТ'')

\noindent \hskip3mm {\bf print}(``Пусть функция $f(z)$ голоморфна в области $W$, определяющейся неравенствами'')

\noindent \hskip3mm {\bf for} $i$ {\bf from} 1 {\bf to} число элементов в $LR$ {\bf do}

\noindent \hskip6mm {\bf print}($W[i]<P[i]$)

\noindent \hskip3mm {\bf end do}

\noindent \hskip3mm {\bf print}(``и непрерывна в замыкании $W$. Тогда в пересечении $D\cap W$, где область $D$, определяется неравенствами'')

\noindent \hskip3mm {\bf for} $i$ {\bf from} 1 {\bf to} число элементов в $z\_s$ {\bf do}

\noindent \hskip6mm {\bf print}($z\_s[i]<Kahler[i]$)

\noindent \hskip3mm {\bf end do}

\noindent \hskip3mm {\bf print}(``справедливо интегральное представление'')

\noindent \hskip3mm {\bf print}($f(\zeta)=\frac1C\int\limits_{\Gamma(\rho)}f(z)\omega(z-\zeta)$)

\noindent \hskip3mm {\bf print}(``где цикл $\Gamma(\rho)$ определяется равенствами'')

\noindent \hskip3mm {\bf for} $i$ {\bf from} 1 {\bf to} число элементов в $LR$ {\bf do}

\noindent \hskip6mm {\bf print}($W[i]=P[i]$)

\noindent \hskip3mm {\bf end do}

\noindent \hskip3mm {\bf print}(``а константа $C=\int\limits_{\Gamma}\omega(z)>0$.'')

\noindent {\bf end procedure}

\medskip

Процедура {\bf LinRel} по заданному набору векторов $vec\_list$ вычисляет все независимые линейные соотношения на них. На выходе процедура выдает список векторов-строк матрицы линейных соотношений.

\medskip

\noindent {\bf LinRel} := {\bf procedure}($vec\_list$)

\noindent \hskip3mm $d$ := число элементов в $vec_list$

\noindent \hskip3mm $n$ := число элементов в $vec_list[1]$

\noindent \hskip3mm {\bf for} $i$ {\bf from} 2 {\bf to} $d$ {\bf do}

\noindent \hskip6mm {\bf if} число элементов в $vec\_list[i]\neq n$ {\bf then} 

\noindent \hskip9mm {\bf print}(``Вектора должны иметь одинаковую длину!'') 

\noindent \hskip9mm {\bf return} 

\noindent \hskip6mm {\bf end if}

\noindent \hskip3mm {\bf end do}

\noindent \hskip3mm $B$ := пустой список

\noindent \hskip3mm {\bf for} $i$ {\bf from} 1 {\bf to} $n$ {\bf do} 

\noindent \hskip6mm добавить в $B$ элемент 0

\noindent \hskip3mm {\bf end do}

\noindent \hskip3mm $b$ := вектор из элементов списка $B$

\noindent \hskip3mm $V$ := матрица, построенная на векторах их $vec\_list$

\noindent \hskip3mm $GenSol$ := общее решение однородной СЛАУ, заданной матрицей $V$ со свободными переменными $s[i]$

\noindent \hskip3mm $S$ := список из всех свободных переменных $GenSol$

\noindent \hskip3mm $FundSol$ := пустое множество

\noindent \hskip3mm {\bf for} $i$ {\bf from} 1 {\bf to} число элементов в $S$ {\bf do}

\noindent \hskip6mm $T[i]$ := 1

\noindent \hskip6mm {\bf for} $j$ {\bf from} 1 {\bf to} число элементов в $S$ - 1 {\bf do}

\noindent \hskip9mm {\bf if} $i+j=$ число элементов в $S$ {\bf then} $k$ := $i+j$ {\bf else} $k$ := $i+j$ по модулю числа элементов в $S$) {\bf end if}

\noindent \hskip9mm $T[k]$ := 0

\noindent \hskip6mm {\bf end do}

\noindent \hskip6mm $a$ := пустой список

\noindent \hskip6mm {\bf for} $j$ {\bf from} 1 {\bf to} число элементов в $S$ {\bf do}

\noindent \hskip9mm добавить в $a$ элемент $T[j]$

\noindent \hskip6mm {\bf end do}

\noindent \hskip6mm $St$ := $S$

\noindent \hskip6mm $k$ := 0

\noindent \hskip6mm {\bf for} $j$ {\bf from} 1 {\bf to} $d$ {\bf do}

\noindent \hskip9mm {\bf if} $St$ зависит от значения $s[j]$ {\bf then} 

\noindent \hskip12mm $St$ := разность множеств $St$ и $s[j]$

\noindent \hskip12mm $k$ := $k+1$

\noindent \hskip12mm $s[j]$ := $a[k]$

\noindent \hskip9mm {\bf end if}

\noindent \hskip6mm {\bf end do}

\noindent \hskip6mm добавить в $FundSol$ значение $GenSol$

\noindent \hskip3mm {\bf end do}

\noindent \hskip3mm {\bf return}($FundSol$, отсортированный по убыванию) 

\noindent {\bf end procedure}

\medskip

Процедура {\bf PrimColl} по заданному набору номеров векторов $cone\_list$, составляющих конусы максимальной размерности в веере, вычисляет все примитивные наборы векторов веера.

\medskip

\noindent {\bf PrimColl} := {\bf procedure}($cone\_list$)

\noindent \hskip3mm $m$ := число элементов в $cone\_list$

\noindent \hskip3mm $n$ := число элементов в $cone\_list[1]$

\noindent \hskip3mm {\bf for} $i$ {\bf from} 2 {\bf to} $m$ {\bf do}

\noindent \hskip6mm {\bf if} число элементов в $cone\_list[i]\neq n$ {\bf then} 

\noindent \hskip9mm {\bf print}(``Все конусы должны иметь одинаковую размерность!'') 

\noindent \hskip9mm {\bf return} 

\noindent \hskip6mm {\bf end if}

\noindent \hskip3mm {\bf end do}

\noindent \hskip3mm $ConeSet$ := множество из элементов $cone\_list$

\noindent \hskip3mm $ConeUnion$ := объединение всех элементов $ConeSet$

\noindent \hskip3mm $AllCones$ := пустое множество

\noindent \hskip3mm {\bf for} $i$ {\bf from} 1 {\bf to} число элементов в $ConeSet$ {\bf do}

\noindent \hskip6mm добавить в $AllCones$ элемент {\bf SetOfAllSubsets}($ConeSet[i]$, 0)

\noindent \hskip3mm {\bf end do}

\noindent \hskip3mm $AllColl$ := {\bf SetOfAllSubsets}($ConeUnion$, 0)

\noindent \hskip3mm $S$ := разность множеств $AllColl$ и $AllCones$

\noindent \hskip3mm $PrimColl$ := пустое множество

\noindent \hskip3mm {\bf for} $i$ {\bf from} 1 {\bf to} число элементов в $S$ {\bf do}

\noindent \hskip6mm {\bf if} {\bf SubsetsMinus1}($S[i]$) является подмножеством $AllCones$ {\bf then} 

\noindent \hskip9mm добавить в $PrimColl$ элемент $S[i]$ 

\noindent \hskip6mm {\bf end if}

\noindent \hskip3mm {\bf end do}

\noindent \hskip3mm {\bf return}($PrimColl$)

\noindent {\bf end procedure}

\medskip

Процедура {\bf nu\_sigma} вычисляет показатели $\nu_l^{\sigma_m}$ по формуле 
\begin{equation}\label{mu}
\nu_l^{\sigma_m}:=-\sum_{i=1}^{n}\det(v_{m_1},\ldots,v_{m_{i-1}},v_l,v_{m_{i+1}},\ldots,v_{m_n}),
\end{equation}
используя процедуру {\bf Permut}, которая, получая на входе вектор и число, выдает список векторов, получающихся из данного путем подстановки заданного числа, на первое, второе, и так далее, места. На входе {\bf nu\_sigma} получает список векторов веера и результат подстановки числа $l$ в вектор $(m_1,\ldots,m_n)$ из (\ref{mu}).

(\ref{mu})\medskip

\noindent {\bf Permut} := {\bf procedure}($c\_list$, $vecnum$)

\noindent \hskip3mm $n$ := число элементов в $c\_list$

\noindent \hskip3mm $NS$ := пустой список

\noindent \hskip3mm {\bf for} $i$ {\bf from} 1 {\bf to} $n$ {\bf do}

\noindent \hskip6mm $t$ := пустой список

\noindent \hskip6mm {\bf for} $j$ {\bf from} 1 {\bf to} $i-1$ {\bf do}

\noindent \hskip9mm добавить в $t$ элемент $c\_list[j]$

\noindent \hskip6mm {\bf end do}

\noindent \hskip6mm добавить в $t$ элемент $vecnum$

\noindent \hskip6mm {\bf for} $j$ {\bf from} $i+1$ {\bf to} $n$ {\bf do}

\noindent \hskip9mm добавить в $t$ элемент $c\_list[j]$

\noindent \hskip6mm {\bf end do}

\noindent \hskip6mm добавить в $NS$ элемент $t$

\noindent \hskip3mm {\bf end do}

\noindent \hskip3mm {\bf return}($NS$)

\noindent {\bf end procedure}

\medskip

\noindent {\bf nu\_sigma} := {\bf procedure}($v\_list$, $l\_list$)

\noindent \hskip3mm $n$ := число элементов в $l\_list$

\noindent \hskip3mm $ns$ := 0

\noindent \hskip3mm {\bf for} $i$ {\bf from} 1 {\bf to} $n$ {\bf do}

\noindent \hskip6mm $t$ := пустой список

\noindent \hskip6mm {\bf for} $j$ {\bf from} 1 {\bf to} $n$ {\bf do}

\noindent \hskip9mm добавить в $t$ элемент $v\_list[l\_list[i][j]]$

\noindent \hskip6mm {\bf end do}

\noindent \hskip6mm $ns$ := $ns$ + определитель матрицы, построенной на векторах из $t$

\noindent \hskip3mm {\bf end do}

\noindent \hskip3mm {\bf return}($(-1)*ns$)

\noindent {\bf end procedure}

\medskip

Процедуры {\bf SubsetsMinus1} и {\bf SetOfAllSubsets} можно отнести к разряду вспомогательных. Процедура {\bf SubsetsMinus1} получает на входе множество и выдает множество всех подмножеств данного множества с числом элементов на единицу меньшим числа элементов входного множества. Процедура {\bf SetOfAllSubsets} получает на входе множество и неотрицательное целое число $K$. Если $K=0$, то процедура выдает множество всех подмножеств данного множества. Если $K>0$, то процедура выдает множество всех подмножеств данного множества, состоящих из $K$ элементов.

\medskip

\noindent {\bf SubsetsMinus1} := {\bf procedure}($inset$)

\noindent \hskip3mm $n$ := число элементов в $inset$

\noindent \hskip3mm $outset$ := пустое множество

\noindent \hskip3mm {\bf for} $i$ {\bf from} 1 {\bf to} $n$ {\bf do}

\noindent \hskip6mm добавить в $outset$ разность множеств $inset$ и $inset[i]$

\noindent \hskip3mm {\bf end do}

\noindent \hskip3mm {\bf return}($outset$)

\noindent {\bf end procedure}

\medskip

\noindent {\bf SetOfAllSubsets} := {\bf procedure}($set$, $K$)

\noindent \hskip3mm $S$ := пустое множество

\noindent \hskip3mm $T$ := множество элементов из $set$

\noindent \hskip3mm {\bf for} $i$ {\bf from} 1 {\bf to} число элементов в $set - K$ {\bf do}

\noindent \hskip6mm $P$ := пустое множество

\noindent \hskip6mm {\bf for} $j$ {\bf from} 1 {\bf to} число элементов в $T$ {\bf do}

\noindent \hskip9mm $R$ := {\bf SubsetsMinus1}($T[j]$)

\noindent \hskip9mm добавить в $P$ элемент $R$

\noindent \hskip6mm {\bf end do}

\noindent \hskip6mm $T$ := $P$

\noindent \hskip6mm добавить в $S$ элемент $T$

\noindent \hskip3mm {\bf end do}

\noindent \hskip3mm {\bf if} $K = 0$ {\bf then} 

\noindent \hskip3mm добавить в $S$ элемент $set$

\noindent \hskip6mm {\bf return}($S$) 

\noindent \hskip6mm {\bf else} {\bf return}($T$) 

\noindent \hskip3mm {\bf end if}

\noindent {\bf end procedure}

\end{document}